\def\strutdepth{\dp\strutbox}
\def \ss{\strut\vadjust{\kern-\strutdepth \sss}}
\def \sss{\vtop to \strutdepth{

\baselineskip\strutdepth\vss\llap{$\diamondsuit\;\;$}\null}}

\newcount\footno
\global\footno=1
\def\comm#1{
\footnote{$^{\the\footno}$}{\sevenrm #1} \advance\footno by 1}

\magnification=\magstep1
\baselineskip16pt
\def\today{\ifcase\month\or
   January\or February\or March\or April\or May\or June\or
    July\or August\or September\or October\or November\or December\fi
\space\number\day, \number\year}
%\line{\hfil\today}
\bigskip
\raggedbottom

\newcount\figno
\global\figno=1
\def\figure#1{\midinsert
\centerline{\epsfbox{#1}}
\centerline{Figure \the\figno}
\global\advance\figno by 1
\endinsert}

\newcount\secno
\newcount\subsecno
\newcount\stno
%\global\stno=1
\global\secno=-1
\outer\def\section#1\par{
     \global\advance\secno by 1
     \global\stno=1
     \global\subsecno=0
    % \vskip0pt plus .3\vsize\penalty-250
     \vskip0pt plus-.3\vsize\bigskip\vskip\parskip
     \message{#1}\noindent{\bf{\S\the\secno.\enspace#1}}
     \nobreak\smallskip\noindent}
\outer\def\subsection#1\par{
      \medbreak
      \global\advance\subsecno by 1
      \noindent{\bf#1}
 %    \global\stno=1
      \nobreak\smallskip}
\outer\def\remark#1\par{
     \medbreak
     \noindent{{\bf Remark   \the\secno.\the\stno.\enspace} {#1}\par
     \global\advance\stno by 1
     \ifdim\lastskip<\medskipamount \removelastskip\penalty55\medskip\fi}}

\outer\def\state#1. #2\par{
     \medbreak
     \noindent{\bf{#1\enspace \the\secno.\the\stno.\enspace}\sl{#2}\par
     \global\advance\stno by 1
     \ifdim\lastskip<\medskipamount \removelastskip\penalty55\medskip\fi}}

\def\tag#1{\edef#1{\the\secno.\the\stno}}
\def\sectag#1{\edef#1{\the\secno.\the\subsecno}}

\newcount\refno
\global\refno=0
\outer\def\ref#1\par{
     \global\advance\refno by 1
     \item{\the\refno.}#1\medskip}

%\font\Bbb=msbm10
%\def\Bbb{\bf}
%%%%%%%%%%%%definitions%%%%%%%%%%%%%%%%%%%%%
\def\QED{\hfill\rlap{$\sqcup$}$\sqcap$\par\bigskip}
\def\give#1.{\medbreak
             \noindent{\bf#1.}}

 %for intersection components
  %for $G$-graphs (changes later to cal G for graphs of groups)
 %for quotient of $\G$-graph by $\G$
 %for maximal invariant forests
 %for invariant subsets of edges
 %for sets of ideal edges in contractibility proof
 % for ideal edges
 % for ideal edges

\def\less{<}

%%%%%%%%%%%%%%%%%%%%%%%%%%%%%%%%%%

\def\{\partial}

\def\int{\text{int}}
\def\invlimit{\smash{\lim\limits_{\raise1pt\hbox{$\longleftarrow$}}}\vphantom{\big(}}
\def\inter{\hskip 1.5pt\raise4pt\hbox{$^\circ$}\kern -1.6ex}
\def\reals{\hbox{I\kern-.13em R}} 
\def\skel(#1,#2){#1^{(#2)}}
\def\hyp {\hbox {\rm {H \kern -2.8ex I}\kern 1.25ex}}
\def\reals {\hbox {\rm {R \kern -2.8ex I}\kern 1.15ex}}
\def\hyp {\hbox {\rm {H \kern -2.8ex I}\kern 1.15ex}}
\def\integers {\hbox {\rm { Z \kern -2.8ex Z}\kern 1.15ex}}
\def\naturals {\hbox {\rm {N \kern -2.8ex I}\kern 1.20ex}}

\def\hyp {\hbox {\rm {H \kern -2.7ex I}\kern 1.25ex}}

\def\HN#1{\parindent0pt \hangafter1 \hangindent35pt \leavevmode%
				 \hbox to 33pt{{\bf #1} \hfil}}
\overfullrule=0pt

\hyphenation{cell-de-com-po-si-tion}
\hyphenation{de-com-po-si-tion}
\hyphenation{di-men-sio-nal}

\font\yoav=cmr8

%input pics.tex

\input BoxedEPS
\SetepsfEPSFSpecial
\HideDisplacementBoxes
%%%%%%%%%%%%%%%%%%%%%%%%%%%%%%%%%%%%%%%%%%%%

\centerline {\bf Additive tunnel number and primitive elements} 
\bigskip

\centerline{ Yoav Moriah\comm{Supported by The Fund for Promoting 
Research at the Technion grant 100-053.}} 
\bigskip

{\yoav{ \noindent {\bf Abstract:} It is proven here that if the connected sum of two tunnel number one 
knots in  $S^3$ is a tunnel number two knot then at least one of the summand knots has a genus two 
Heegaard splitting with a meridian as a primitive element. Hence this is a necessary and sufficient 
condition for tunnel  number one knots to have additive tunnel number.}}
\bigskip
\bigskip
\section Introduction

\bigskip
\noindent
The way in which the tunnel number $t(K)$ of a knot $K = K_1 \# K_2$   relates to $t(K_1)$ and 
$t(K_2)$  is a long standing question. It was long known that $t(K_1 \# K_2) \leq t(K_1) + t(K_2) +1$. 
That this inequality is best possible was proved by the author and Rubinstein in  [MR] and 
by  Morimoto, Sakuma and Yokota in [MSY]. In the other direction it was  shown by Morimoto in 
[Mo1] that there are prime knots for which $t(K_1 \# K_2) \less t(K_1) + t(K_2) $. For tunnel number 
one knots the sum cannot decrease since it was shown by Norwood [No] that tunnel number one knots 
are prime. It is also known (see Section 1) that if one of the two knots has a genus two Heegaard 
splitting in which a meridian is primitive then the tunnel number is additive. It remained an open
question whether this is also a necessary condition as it is conceivable that the exterior of the connected
sum of the knots has a genus two Heegaard splitting which is not induced by any Heegaard splitting
of the summand knots.

\bigskip

\noindent The main theorem of this paper is the following:
\bigskip

\tag\mainThm

\state Theorem. Let $K_1$ and $K_2$ be  tunnel number one knots  in $S^3$ . Assume that  $t(K) =
 t(K_1) + t(K_2)$. Then at least one of $K_1$ or $K_2$ has a genus two Heegaard splitting in which 
a meridian curve represents a primitive element in the handlebody component of the splitting.

\bigskip

As an immediate corollary we obtain:
\bigskip

\tag\necessaryCor

\state Corollary.  Let $K_1$ and $K_2$ be  tunnel number one knots  in $S^3$. Then $t(K) =
 t(K_1) + t(K_2)$ if and only if one of $K_1$ or $K_2$ has a genus two Heegaard splitting in which 
a meridian curve represents a primitive element in the handlebody component of the splitting.
\bigskip

\eject

\noindent

 For definitions of the above terminology see Section 1.

\bigskip

\noindent {\bf Acknowledgements:}  I would like to thank Cameron Gordon  and John Luecke for 
conversations regarding this work and the University of Texas at Austin, where this
work was done, for its hospitality. 
\bigskip
\bigskip

\section Preliminaries

\bigskip
\noindent In this section we define some of the notions and prove some technical lemmas needed
for the proof of the main theorem.

Throughout the paper  $K_1$ and $K_2$ will be knots in $S^3$ and  $K = K_1 \# K_2$ will denote 
the connected sum of $K_1$ and $K_2$ (unless specified otherwise).  Let $N()$ denote an open 
regular neighborhood in $S^3$.

Recall that $(S^3,K)$ is obtained by removing  from each space $(S^3,K_i),  i = {1,2}$ a 
small $3$-ball intersecting $K_i$ in a short unknotted arc and gluing the two remaining $3$-balls along 
the $2$-sphere boundary so that the pair of points of $K_1$ on the $2$-sphere are identified with  the 
pair of points of $K_2$. If we denote $S^3 - N(K)$ by $E(K)$ then $E(K)$ is obtained from 
$E(K_i), i = {1,2}$ by identifying a meridional annulus $A_1$ on $\partial E(K_1)$ with a 
meridional annulus $A_2$ on $\partial E(K_2)$. A knot $K \subset S^3$ is {\it prime} if it not a 
connected sum of non-trivial knots. The annulus $A_1 = A_2$ will be denoted by $A$ 
and called the{ \it decomposing annulus}.

 A {\it tunnel system} for a arbitrary knot $K \subset S^3$ is a collection of 
properly  embedded locally unknotted arcs  $t_1, \dots, t_n$  in  $S^3 - N(K)$  so that 
$S^3 - N(K \cup t_1 \cup \dots \cup t_n)$ is a handlebody.

Given a tunnel system for a knot $K \subset S^3$ note that the closure of 
$N(K \cup t_1 \cup \dots \cup t_n)$ is always a handlebody denoted by $V_1$ 
and the handlebody $S^3 - N(K \cup t_1 \cup  \dots \cup t_n)$ will be denoted by $V_2$. 
For a given knot $K \subset S^3$ the smallest cardinallity of any tunnel system is called the 
{\it tunnel number} of $K$ and
is denoted by $t(K)$.

A compression body  $V$  is a  3-manifold with a preferred boundary component
$\partial_+V$  and is obtained from a collar of  $\partial_+ V$  by attaching 2-handles and
3-handles, so that the connected components of  $\partial_- V$ = $\partial V - \partial_+ V$ are all
distinct from  $S^2$.  The extreme cases, where  $V$  is a handlebody i.e., $\partial_- V = \emptyset$,
or where $V = \partial_+V \times I$, are admitted.  Alternatively we can think of $V$ as obtained from
$(\partial_-V) \times I$ by attaching $1$-handles to $(\partial_-V) \times \{1\}$. An annulus in a 
compression body will be called a {\it vertical annulus} if it has its boundary components on different
boundary components of the compression body.

Given a knot  $K \subset S^3$ a {\it Heegaard splitting } for $E(K)$ is a decomposition of $E(K)$ into
a compression body $V_1$ containing $\partial E(K)$ and a handlebody $S^3 - int(V_1)$. Hence, a 
tunnel system $t_1, \dots, t_n$  in  $S^3 - N(K)$  for $K$ determines a Heegaard splitting of genus
$n +1$ for $E(K)$. 

Given a Heegaard splitting $(V_1, V_2)$  for $ S^3 - N(K_1 \# K_2)$ we can assume that the 
decomposing annulus $A$  intersects the compression body $V_1$ in two vertical annuli  
$A^*_1, A^*_2$ and a collection of disks $D_1, \dots, D_l$ . Note also that $A$ intersects 
$V_2$ in a connected planar surface.

Let $ {\cal E} = \{E_1 , \dots, E_{t(K) + 1}\}$ be a complete meridian disk 
system for $V_2$, chosen to minimize the intersection ${\cal E}\cap A$. Since $V_2$ is a handlebody 
it is irreducible  and we can assume that no component of  ${\cal E} \cap A$ is a simple closed curve. 
Furthermore each arc of intersection $\alpha$ in ${\cal E} \cap A$ is an essential arc.

When we cut $E(K)$ along the decomposing annulus $A$ any Heegaard splitting $(V_1, V_2)$
induces  Heegaard splittings on both of $E(K_1)$ and $E(K_2)$: Set $V_1^i = (V_1 \cap E(K_i))
\cup N(A)$ and $V_2^i = V_2 - N(A)$. The pair $(V_1^i,V_2^i)$ is a Heegaard splitting for $E(K_i)$.

We say that an element $x $  in a free group $F_n$ is {\it primitive} if it belongs to some basis for 
$F_n$. A curve on a handlebody $H$  is {\it primitive} if it represents a primitive element in the free
group $\pi_1(H)$. An annulus $A$ on  $H$ is {\it primitive} if its core curve is primitive. Note that
a curve on a handlebody is primitive if and only if there is an essential disk in the handlebody 
intersecting the curve in a single point.

Two Heegaard splittings $(V_1^i, V_2^i)$ for  $E(K_i)$ repectfully, induce a decomposition of
$E(K)$ into $(V_1, V_2)$. We obtain $V_1$ by gluing the compression bodies $V_1^1$ and 
$V_1^2$ along two vertical annuli and $V_2$ by gluing $V_2^1$ and $V_2^2$ along a meridional
annulus. Hence $V_1$ is always a compression body but $V_2$ is a handlebody if and only if
the meridional annulus is a primitive annulus.

Following some ideas of Morimoto (see [Mo 2]) we consider now the planar surface $ P = A \cap V_2$.
It has two distinguished boundary components coming from the vertical annuli  $A^*_1, A^*_2$ and 
denoted by $C^*_1, C^*_2$ respectively. There are exactly $d$ other boundary components of $P$ which we denote 
by $C_1, \dots, C_d$. With this notation we have $\partial D_i = C_i$.  The arcs of $E \cap A$ are 
contained in $P$ and come in three types: 
{\item 1.  An arc $\alpha$ of Type I is an arc connecting two different boundary components of $P$.}
{\item 2.  An arc $\alpha$ of Type II is an arc connecting a single boundary component of $P$ 
to itself 

\hskip8pt so that the arc does not separate the boundary components $C^*_1, C^*_2$.}

{\item 3.  An arc $\alpha$ of Type III is an arc connecting a single boundary component of $P$ 
to 

\hskip8pt itself. with the additional property that the arc does separate the boundary com-

\hskip8pt ponents $C^*_1, C^*_2$.}
\smallskip
Since the $A$ was chosen to minimize the number of disks in $A \cap V_1$ then $P$ is incompressible
in the handlebody $V_2$. Hence there is a sequence of boundary compressions of $P$ along 
disjoint arcs $\alpha_i$ using sub-disks of $\cal E$ so that the end result is a collection of disks. 
Any such sequence defines an order on the  arcs $\alpha_i$ .

\smallskip

\tag\darcDef
\smallskip
\state Definition. We call $\alpha_i$, an arc of intersection  of $P \cap E$, a  $ d-arc$ if $\alpha_i$ is 
of type I and there is some component $C$ of $\partial P - (C^*_1\cup C^*_2)$  which  meets $\alpha_i$
and does not meet any $\alpha_j$ for any $j \less i$. If  $\alpha_i$ is of type I and connects $C^*_1$ to
$C^*_2$ it is called an $e-arc$

Any outermost arc $\alpha_i$ determines a sub-disk $\Delta$ on some $E_i$ where $\partial \Delta =
\alpha_i \cup \beta$ and $\beta$ is an arc on $\partial V_1 = \partial V_2$. When we perform an 
isotopy of type A  i.e., pushing $P$ through $\Delta$ as in [Ja], we produce a band $b$ with core
$\beta$ on $\partial V_1 = \partial V_2$. The following crucial result is proved in  [Mo2] pp 41- 42, 
and [Oc]: 
\smallskip

\tag\darcThm
\bigskip
 \state Theorem.  (Morimoto) If the decomposing annulus is chosen to minimize the number of
components of $V_1 \cap A$ and $V_1 \cap A \not= \emptyset$  then in $V_2 \cap E$ = $P \cap E$:
\item {\rm (a)} there are no d-arcs. 
\item {\rm (b)} there are no e-arcs.
\item {\rm (c)} there are no arcs of type II.
\item {\rm (d)} each component $C \subset \partial P$ has an arc $\alpha$ of type III with end points
on $C$. 

\bigskip
Consider now a Heegaard splitting $(V_1 , V_2)$  for $E(K)$ the exterior of $K = K_1 \# K_2$, 
where  $\partial E(K) \subset V_1$ and  the decomposing annulus $A$ meets $V_1$ in disks 
and two vertical annuli.  Since the annulus $A$  meets $V_2$ in a connected planar surface $P$
it separates $V_2$ into two components each of which is a handlebody. We will denote the 
handlebodies $cl(V_2 - A) \cap E(K_i)$ by $V_2^i$ respectively. 
However $V_1 - A $ might have many components. 

\bigskip
\tag\floatDef
\bigskip
\state Definition. A component of $cl(V_1 - A) $ which is disjoint from $\partial E(K_i)$ and 
intersects $A$ in $n$ disks will be called a  {\it n-float}.

\bigskip

\noindent {\bf Remark:} Note that a n-float is either a 3-ball or a handlebody if its spine is 
not a tree. Furthermore there is always exactly two components of $cl(V_1 - A) $ not disjoint 
from $\partial E(K_i)$ (one in each of $E(K_1)$ and $E(K_2)$) and each one is a handlebody 
of genus at least one as $V_1$ is a compression body with a $T^2$ boundary. We denote these 
special  components by $N_1$ and $N_2$ depending on whether they are contained in 
$E(K_1)$ or $E(K_2)$ respectively.

\bigskip
Consider now $E_i \subset \cal E$ any one of the meridian disks of $V_2$. On $E_i$ we have a 
collection of arcs corresponding to the intersection with the decomposing annulus. These arcs, 
as indicated in Fig. 1 below, separate $E_i$ into sub-disks where disks on opposite sides of arcs  
are contained in opposite sides of $A$ i.e., in $E(K_1)$ or  $E(K_2)$ respectively. So each sub-disk 
is contained in either $E(K_1)$ or  $E(K_2)$. The boundary of these sub-disks is a collection 
of alternating arcs $\cup(\alpha_i  \cup \beta_i)$ where $\alpha_i$ are arcs on $A$ and $\beta_i$
are arcs on some component of  $cl(V_1 - A)$.
\bigskip

\bigskip
\hskip115pt \BoxedEPSF{Fig1.eps scaled 500}
%\hskip105pt { \picture 2.00in by 1.80in  (Fig.1 scaled 500)}
\bigskip
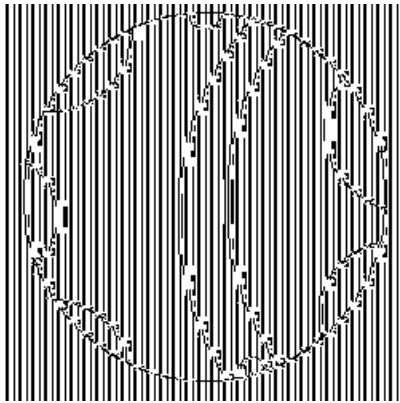
\centerline {Fig. 1}
\bigskip

\bigskip
\tag\innertunnelProp
\bigskip
  \state Proposition. Let $K_1$ and $K_2$ be knots in $S^3$ and let $K, A , \cal E$ be the  connected
sum, minimal intersection decomposing annulus and meridional system for some Heegaard splitting
of $E(K)$ as above. Then 
\item {\rm (a)} the $\beta$ arc part of the boundary of an outermost sub-disk in $E$ cannot
 be contained in a n-float which has no genus.
\item {\rm (b)}  if  the $\beta$ arc part of the boundary of an outermost 
sub-disk in $E$ is contained in a $N_i$ component $i = 1$ or $2$ the genus of $N_i$ is greater than 
one.
\bigskip

\give Proof. Denote an outermost sub-disk of some $E_j$ by $ \Delta$  and suppose it is cut off
by an arc $\alpha$ on $A$ with end points of a disk $D_i$  which belongs to some n-float which
has no genus. Further assume $\partial \Delta = \alpha \cup \beta$ where $\beta$ is an arc on the 
n-float meeting $D_i$ in exactly two points $\partial \beta = \partial \alpha$. On $\partial D_i$ 
there is a small arc $\gamma$ so that $ \gamma \cup \beta$ is a simple closed curve 
on the n-float bounding a disk $D$ there, since the n-float has no genus (see Fig. 2 below).
Furthermore $\gamma \cup \alpha$ is a simple closed loop on $A$ which bounds a sub-annulus 
of $A$. Hence $\gamma \cup \alpha$ bounds a disk $D'$ on the decomposing $2$-sphere of $K$
intersecting $K$ in a single point. Thus we obtain a $2$-sphere $D \cup \Delta \cup D'$ which 
intersects the knot $K$ in a single point. This is a contradiction  which finishes case (a). 

For case (b), assume that  the outermost disk $\Delta$ is contained in $N_1$, say, and that genus 
$N_1$ is one. As before we have $\partial \Delta = \alpha \cup \beta$ where $\beta$ is an arc on 
$N_1$ and a small arc $\gamma$ so that $ \gamma \cup \beta$ is a simple closed curve on $N_1$. 
If $ \gamma \cup \beta$ bounds a disk in $N_1$ we have the same proof as in case (a). If
 $ \gamma \cup \beta$ does not bound a disk on $N_1$ we consider small sub-arcs $\beta_1$
and $\beta_2$ of $\beta$ which are respective closed neighborhoods of $\partial \beta$. These
arcs together with a small arc $\delta$ on $\partial N_1  - \partial E(K_1)$  and $\gamma$ bound 
a small band $b$ on $\partial N_1$. Notice that $b  \cup_ {\beta_1,\beta_2}\Delta$ is an annulus
$A'$. The annulus $A'$ together with the sub-annulus $A''$ of $A$ cut off by $\alpha \cup \gamma$ 
defines an annulus $A' \cup_{\alpha \cup \gamma} A''$ which determines an isotopy of a meridian 
curve $C_1$ to a simple closed curve $\lambda$ on $\partial N_1$. Note that $N_1$ is a solid torus 
and $\pi_1(N_1) =  \integers$ which is generated by a meridian $\mu$ of $E(K_1)$. Hence 
$[\lambda] =   [C_1] = \mu \in \pi_1(N_1)$. (see Fig. 3). Now we can consider the  
annulus $(A - A'') \cup A'$.  If it is non-boundary parallel it is a  decomposing annulus with at least 
one less disk component intersection than $A$ in contradiction to the choice of $A$. If it is boundary
parallel we have $A'' \cup A'$ as a decomposing annulus with a smaller  number of disks. Again in contradiction 
to the choice of $A$. So genus $N_1$ cannot be one and this finishes case (b).

\QED 
\hskip50pt \BoxedEPSF{Fig2.eps scaled 550}
%\hskip50pt  { \picture 4.00in by 2.80in  (Fig.2 scaled 550)}
\bigskip
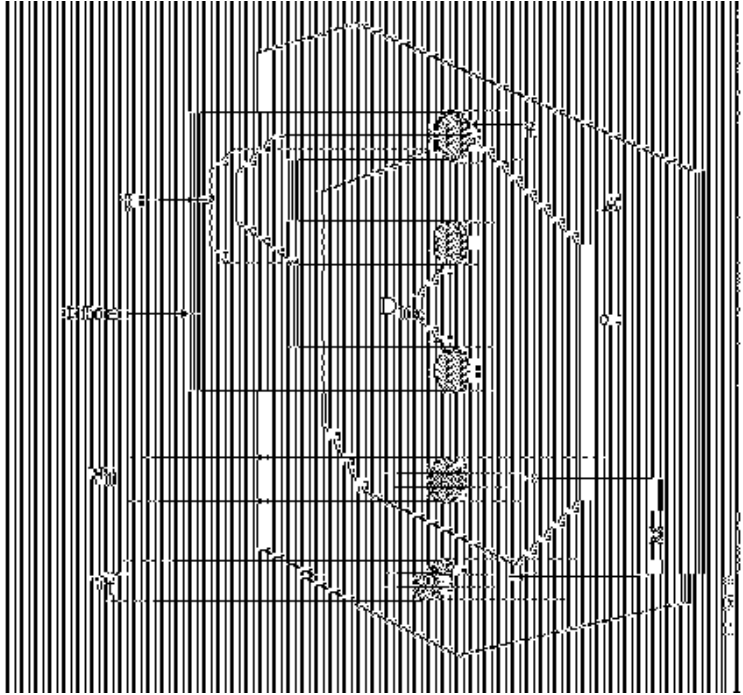
\centerline {Fig. 2}
\bigskip
\bigskip
\bigskip
\hskip50pt \BoxedEPSF{Fig3.eps scaled 550}
%\hskip50pt { \picture 4.00in by 2.80in  (Fig.3 scaled 550)}
\bigskip
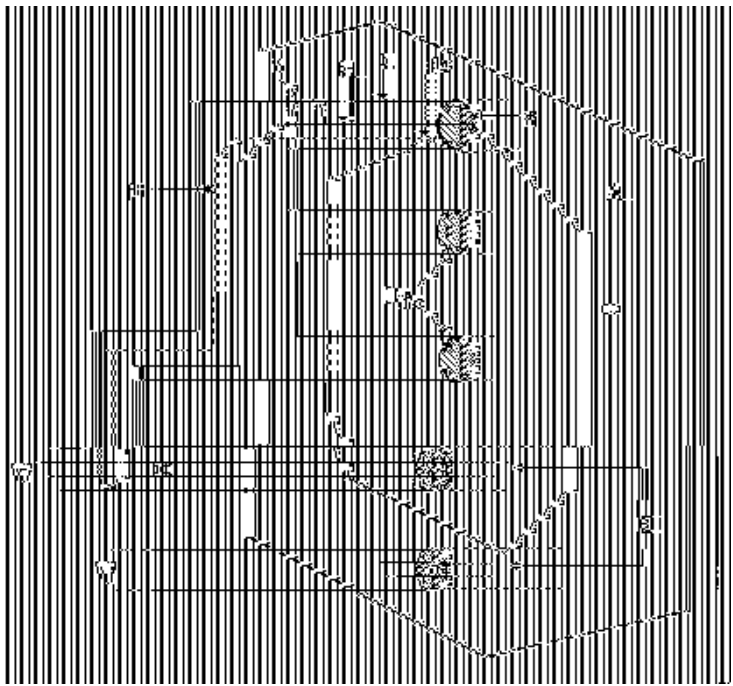
\centerline {Fig. 3}
\bigskip

\tag\innertunnelCor
\bigskip
  \state Corollary. Every unknotting tunnel system $\tau$ for $K= K_1 \# K_2$ must contain at least 
one tunnel which is disjoint from a decomposing annulus for $K$ minimizing the number of 
intersections with  $N(K\cup \tau)$.
\QED

We finish this section with a proposition which is probably due to Morimoto and Sakuma [MS] 
as it is implicit in their discussion of dual tunnels there (see also [Sc]). 

\smallskip
\tag\tunnelrepProp
\smallskip
\state Proposition. Let  $(V_1, V_2)$ be a Heegaard splitting of a $3$-manifold $M$. Assume 
that there is an incompressible annulus $A \subset V_2$ so that  both boundary components of
$A$ are in $\partial V_2 =  \partial V_1$  and one of the boundary 
components intersects an essential disk $D$ of $V_1$. Then if we add a one handle to $V_1$ 
which is a regular neighborhood of an essential arc in $A$ and remove a regular neighborhood 
of $D$ from $V_1$ we obtain new  compression bodies $V'_1$ and $V'_2 = M- int(V'_1)$
defining a new Heegaard splitting of the same genus of $M$.

\smallskip

\give Proof. It is clear that $V'_1$ is a compression body. To see that $V'_2$ is a compression body
we go through an intermediate step. Cutting $V_2$ along $A$ we obtain a compression body as
$A$ is incompressible and since $\partial A$ meets $D$ is a single point $A$ is non-separating.
If we add a regular neighborhood of an essential arc $ \nu$ of $A$  to $V_1$ we obtain a compression 
body $V^\#_1$ of  genus one bigger than that of $V_1$. Its complement is obtained by gluing together
the two copies of the disk $A - N(\nu)$ on the compression body. Thus getting a compression body
$V^\#_2$ of genus one bigger than that of $V_2$. However in the Heegaard splitting 
$(V^\#_1, V^\#_2)$ the essential disk $A - N(\nu)$ intersects the essential disk $D$ in a single point. So
we can reduce the genus by one by removing a regular neighborhood of $D$ from $V^\#_1$ and adding
 it to $V^\#_2$. Thus obtaining the compression bodies $V'_1$ and $V'_2$ (as indicated in Fig. 4).
 
\QED

\vskip80pt
\hskip40pt \BoxedEPSF{Fig4.eps scaled500}

%\hskip25pt { \picture 4.70in by 1.70in  (Fig.4 scaled 500)}

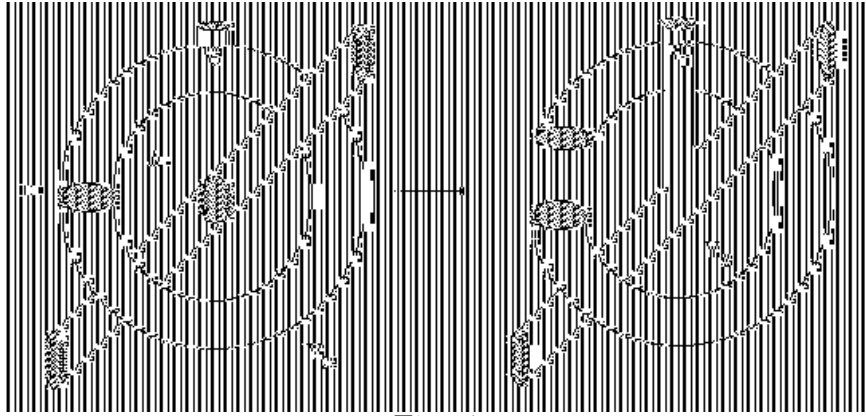
\centerline {Fig. 4}
\bigskip
\bigskip

\section Ruling out U-turns
 
\bigskip
\bigskip

We have the following proposition due to 
Anna Klebanov (see [Kl]). For the completeness of the argument we present it here with proof.

\bigskip
\tag\UturnProp
\bigskip
\state Proposition. (Klebanov) Given a decomposing annulus $A$ which minimizes the number of disk 
components of $V_1 \cap A $ then no component of $cl(V_1 - A) $ is a $3$-ball meeting $A$ 
in exactly two disks.

\noindent {\bf Remark:} In other words, with the above assumptions no tunnel of $K$ has U-turns.
I.e.  no tunnel  pierces $A$ in one direction and then turns around and pierces it again in the opposite
direction without meeting any other part of $V_1$. 

\bigskip

\give Proof. Assume in contradiction that some 2-float is a $3$-ball meeting $A$ 
in exactly two disks. Denote these disks by $D_1$ and $D_2$ but note that these indices do not
necessarily agree with the natural order defined on the disks $D_i$ by  Theorem \darcThm\/.  We
first need the following lemma:

\bigskip
\tag\DiskLemma
\state Lemma. There is some $E_j \subset {\cal E}$ and a sub-disk $\Delta \subset E_j$ so that
$\partial \Delta = \cup (\alpha_r \cup \beta_s)$, where the $\beta_s$ arcs are contained in
the 2-float and the $\alpha_r $ are arcs on $A$  of type I or  III. Hence $\partial \Delta \subset 2-float 
\cup A$.

\give Proof of Lemma. By  Theorem \darcThm\/ there is some arc $\alpha$ with end points on 
$D_1$. This arc $\alpha$ occurs in some $E_j$ and separates it into two disks. Consider the 
disk adjacent to $\alpha$ on the same side of $A$ as the 2-float. We can assume that the sub-disk 
is to the right of $\alpha$. This sub-disk  cannot be outer-most as the 2-float has no genus.  
Hence there are more arcs of intersection further to the right of $\alpha$. If all arcs on the 
sub-disk $\Delta$ adjacent to $\alpha$, which are not on $\partial E_j$, are of type III the disk 
$\Delta$ satisfies the conclusion of the lemma and we are done: Since all arcs of type III have end
points on $D_1$ or $D_2$ the $\beta$ arcs must be contained in the 2-float.

So we assume that some arc to the right of $\alpha$ is of type I. If all such arcs have both end 
points on $D_1$ or $D_2$ we are done as well. The only way in which $\partial \Delta$ can 
leave the 2-float is by means of an arc of type I. Hence if the lemma fails there are at least two 
arcs of type I  
with exactly one end point not on $D_1$ or $D_2$. Since $\partial \Delta$ is connected and if we 
leave the 2-float by means of an arc of type I we must come back to it also by means of an arc of 
type I. Consider such an arc $\rho$, it cannot be outer-most as it is of type I and would  then 
be a  d-arc contradicting Theorem \darcThm\/. Hence further to the right there is some arc $\alpha'$
of type III with end  points on $D_1$ or $D_2$. One of the two disks adjacent to $\alpha'$ is on
the 2-float. Assume that it is $\Delta'$ and that it is on the right of $\alpha'$. Then we start
the argument again with $\alpha'$. This procedure must end since the intersection is finite. 

Assume therefore that $\Delta'$ is to the left of $\alpha'$. If all arcs in $\partial \Delta' -\partial E_j$
are either all of type I or  type III with end points on $D_1$ or $D_2$ we are done as before.  
If there an arc of type I with no end points on $D_1$ or $D_2$ then there is an arc $\rho'$ of type I 
with exactly one end point on $D_1$ or $D_2$.
It cannot be outermost as before so farther out there is an arc of type III with end points on
$D_1$ or $D_2$. So we can start the argument with  $\rho'$. However the procedure must
terminate as the intersection is finite. Hence at some stage we obtain a disk $\Delta$ with
$\partial \Delta - \partial E_j$ consisting of arcs of type I or type III all of which have 
end points on $D_1$ or $D_2$. Hence all the $\beta$ arcs are on the 2-float. (see Fig.1)

\QED

Consider now an essential sub-annulus $A'$ of $A$ containing the disks $D_1$ and 
$D_2$. It is a meridional annulus in $(S^3, K)$ so we can cap off $A'$ by two meridian disks 
$D^*_1$ and $D^*_2$ in $(S^3, K)$ to obtain a 2-sphere intersecting $K$ in exactly two points in
$D^*_1$ and  $D^*_2$. If we attach the boundary of the 2-float to this 2-sphere along 
$D_1$ and $D_2$ we get a 2-torus $T$. By the above lemma $\partial \Delta$ is contained in $T$.

\bigskip
\tag\essarcLemma
\state Lemma. The loop $\partial \Delta$ is essential in $T$.
\bigskip

\give Proof of Lemma. Assume that $\partial \Delta$ bounds a disk $\Delta'$ on $T$. If  $\Delta'$
contains only one of $D^*_1$ or $D^*_2$ then the 2-sphere $\Delta' \cup \Delta$ intersects $K$
in a single point. If  $\Delta'$ contains both of $D^*_1$ and $D^*_2$ then the 2-sphere 
$\Delta' \cup \Delta$ is a decomposing 2-sphere for $K$ as if it was boundary parallel we could
isotope $A$ off $\Delta$ reducing the intersection of $A$ and $\cal E$ incontradiction.
The 2-sphere $\Delta' \cup \Delta$ meets $V_1$ in less disks than $A$ as  $ \Delta$ is in $V_2$ 
and  $\Delta'$ does not contain  $D_1$ or $D_2$ in contradiction to the choice of $A$.

Since $\partial \Delta$ bounds a disk $\Delta'$ on $T$ the intersection of $\partial \Delta$ with
a core curve of the meridional annulus $A'$ is even. Similarly the intersection of $\partial \Delta$ 
with the boundary of a  cocore disk of the 2-float is even. Hence the number of arcs of type I is
even and so is the number of $\beta$ arcs (these are the arcs which intersect the boundary of a  
cocore disk of the 2-float). Hence the number of arcs of type III (the $\alpha$ arcs) is also even.
As a consequence the disk   $\Delta'$ is a union of bands glued together to each other at their ends.
The bands correspond to the areas in $A$ between the arcs of type I and between the arcs of type III
and also on the 2-float between the $\beta$ arcs (see Fig. 5).

\bigskip
\vskip5pt

\hskip10pt\BoxedEPSF{Fig5.eps scaled 550}
%\hskip0pt { \picture 5.00in by 1.70in  (Fig.5 scaled 590)}

\smallskip
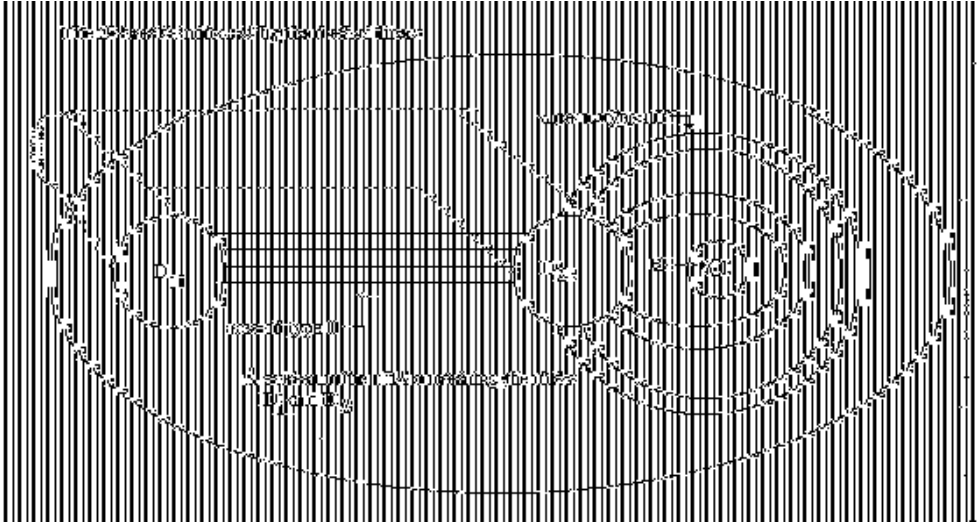
\centerline {Fig. 5} 
\bigskip
Since the bands are glued to each other along small arcs on both ends, the number of gluing arcs
is equal to the number of bands. So an Euler characteristic argument shows that

\smallskip
\centerline {$ \chi (\Delta') = \sum \chi(bands) -  \sum \chi (gluing \hskip3pt arcs) = 0$} 

\bigskip

\noindent  But this is obviously a contradiction and hence $\partial \Delta$ is essential in $T$.

\QED

\bigskip
Now we do 2-surgery on $T$ along the curve $\partial \Delta$ by removing an annulus 
neighborhood of $\partial \Delta$ on $T$ and gluing two copies of $ \Delta$. By an
Euler characteristic argument we obtain a 2-sphere intersecting $K$ in two points
on  $D^*_1$ and $D^*_2$. If we remove  $D^*_1$ and $D^*_2$ we obtain an annulus
$A''$. We can now replace the annulus $A'$ by the annulus $A''$ and get a new decomposing
annulus $(A - A') \cup A''$ which does not intersect the disks $D_1$ and $D_2$ (since $A''$
does not). This contradicts the choice of $A$ and hence we cannot have 2-floats of genus zero
as stated in the proposition. \QED

\bigskip

\section  Additive tunnel number one knots

\bigskip

In this section we assume that both $K_1$ and 
$K_2$ are tunnel number one knots and hence are both prime by a result of Norwood (see [No]). 
We further assume that $t(K) = t(K_1) + t(K_2) = 2$. We have the following lemma:
 \smallskip
\tag\configLemma
\state Lemma. Let $K_1$ and  $K_2$ be tunnel number one knots in $S^3$ so that $t(K) = 
t(K_1) + t(K_2) = 2$. Let $(V_1,V_2)$ be a genus three Heegaard splitting of $S^3 - N(K)$.
If $A$ is a decomposing annulus minimizing intersection with $V_1$ then we can choose a 
spine for $V_1$ which is $\partial E(K) \cup t_1\cup t_2$ where $t_1 \cap A = \emptyset$
and $N(t_2) \cap A$ is either empty or is composed of at most two disks. As indicated in
Fig. 6 (a), (b), (c) and (d).

\bigskip
\vskip10pt

\hskip50pt  \BoxedEPSF{Fig6.eps scaled 590}

%\hskip50pt { \picture 2.40in by 2.40in  (Fig.6 scaled 590)}
\vskip3pt
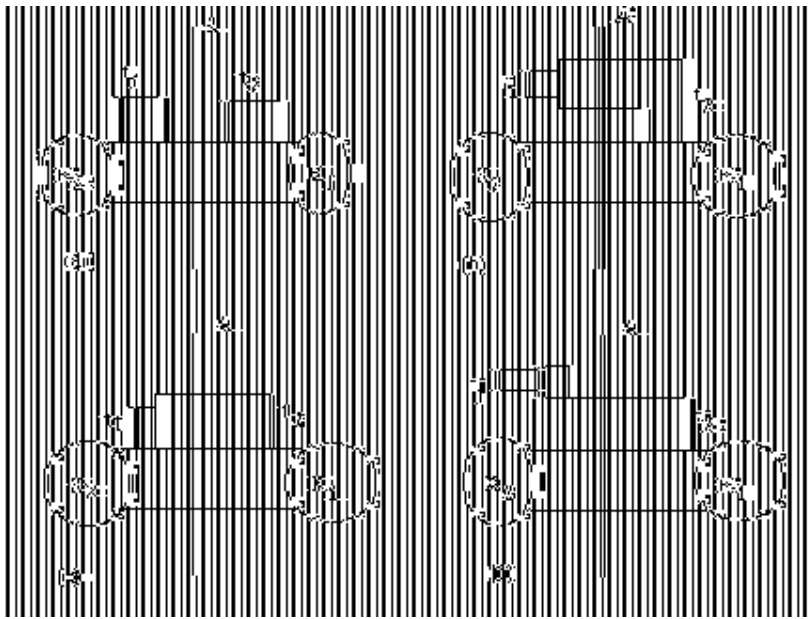
\centerline {Fig. 6}

\eject
\bigskip
\give Proof. \rm By Proposition \innertunnelProp \/ there is a tunnel $t_1$ i.e., a  maximal cocore arc
of an essential disk in $V_1$ which does not meet $A$ at all. If  we cut $V_1$ along the essential disk
we have a  genus two compression body.  We can now choose a spine for the $1$-handle connected
to the $T^2 \times I$ part of the compression body. Denote this arc by $t_2$.  It cannot intersect $A$ in
more than two points as this would create  $2$-floats  of genus  zero  in contradiction to 
Proposition \UturnProp \/. Hence there are four possibilities. 

\smallskip

\noindent 1. The arc $t_2$ has both its end points on one side of $A$ and does not meet $A$.
In this case $t_1$ is on the other side of $A$. Since both knots are tunnel number one knots.

\smallskip

\noindent 2. The arc $t_2$ has both its end points on one side of $A$ and does meet $A$. In this
case $t_2$ meets $A$ in two points and $t_1$ has both of its end points on $t_2$. As otherwise $t_2$
would create a $2$-float of genus zero or if $t_1$ had one end point on $t_2$ and the other on $K$
then $V_1 - A$ would have exactly two components both of genus one in contradiction to Proposition 
\innertunnelProp \/.

\smallskip
\noindent 3. The arc $t_2$ has one end point on one side of $A$ and the other on the other side. 
In this case $t_2$ meets $A$ in a single point as otherwise we have genus zero  $2$-floats. 
The arc $t_1$ can either have both its end points on $t_2$ or one  end point on $t_2$ and 
the other on $K$.

\smallskip
 
\noindent 4. The arc $t_2$ has one end point on one side of $A$ and the other on the other side on
 $t_2$ creating a little loop. In this case $t_1$ must have both its end points on $t_2$ or otherwise
we are in case 3. 

\smallskip

\noindent  All these cases are indicated in Fig. 6

\QED

We are now ready to prove the theorem.

\bigskip

\give Proof.  By Lemma \configLemma\/ we only need to consider the four possible configurations
of Cases $1 - 4$  discussed there. 
\bigskip
\noindent $ \underline{Case  \hskip3pt 1}$:   Since $A$ is incompressible it cuts  $V_2$ into two 
handlebodies $V^1_2$ and  $V^2_2$.  These handlebodies when glued along $A$ yield 
a handlebody.  Hence $A$ must be primitive in one of them.   

\bigskip
\noindent  $\underline{Case \hskip3pt 2}$:  
When we cut $V_2$ along $A$ we obtain a $2$-float of genus one intersecting the annulus
$A$ in two disks denoted by $D_1$ and $D_2$. We can assume that the $2$-float, which we denote
by $V$ is contained in $E(K_2)$. Any outermost disk $\Delta$ in any meridian disk of $V_2$  must 
have $\partial \Delta = \alpha \cup \beta$ where $\alpha$ is an arc of type III and $\beta$ is an arc
on $V$. The shared end points of  $\alpha$  and $\beta$ are either on $D_1$ or $D_2$ and we can
assume that they are on $D_1$ as the picture is symmetric (see Fig. 7, below).
\bigskip
\bigskip
\hskip40pt \BoxedEPSF{Fig7.eps scaled 550}
%\hskip40pt { \picture 4.00in by 2.80in  (Fig.7 scaled 550)}
\smallskip
\centerline {Fig. 7}

\bigskip
Hence we have a solid torus $V$ with two marked disks $D_1$ and $D_2$ on it and there is an
essential arc $\beta$ connecting $D_1$ to itself (see Fig. 8).
\bigskip

\hskip80pt \BoxedEPSF{Fig8.eps scaled 600}

%\hskip70pt { \picture 2.20in by 2.20in  (Fig.8 scaled 600)}
\bigskip

\centerline {Fig. 8}

\vfill \eject
\bigskip
\noindent  $\underline {Claim}$:   There is a meridian disk $D$ of $V$ which intersects 
$\beta \cup D_1$ in a single point.
\bigskip
\noindent {\it Proof  of Claim:} Consider a small arc $\gamma$ on $\partial D_1$ so that
$\gamma \cup \alpha$ is a meridian curve of $E(K_2)$, chosen so that $\gamma \cup \alpha$ 
does  not separate  the disk $D_1$ from $D_2$. As before, consider small sub-arcs $\beta_1$
and $\beta_2$ of $\beta$ which are respective closed neighborhoods of $\partial \beta$. These
arcs together with a small arc $\delta$ on $\partial V  - (\partial D_1 \cup \partial D_2)$  and $\gamma$ 
bound a small band $b$ on $\partial V$. Notice that $b  \cup_ {\beta_1,\beta_2}\Delta$ is an annulus
$A'$. The annulus $A'$ together with the sub-annulus $A''$ of $A$ cut off by $\alpha \cup \gamma$ 
define an annulus $A' \cup_{\alpha \cup \gamma} A''$ which determines an isotopy of a meridian 
curve $C_1 \subset \partial A$ to a simple closed curve $\lambda = \beta - (\beta_1\cup \beta_2) \cup 
\delta$ on $\partial V$.  So $\lambda$ is homotopic in $V$ to some power of a core curve $v$ of $V$
(i.e., $\lambda = v^n, n \in \integers$). Hence in $E(K_2)$ a meridian curve 
$\mu = \alpha \cup \gamma \simeq v^n$. Since $[\mu]$ is a generator of $H_1(E(K_2)) = \integers$
this can only happen if $n = \pm 1$. Hence $\lambda$ is homotopic to a curve on $\partial V$ 
which intersects a meridian disk of $V$ in a single point. As $\partial V = T^2$, it follows that 
$\lambda$ is isotopic to a longitude curve of $V$ and hence there is a disk isotopic to a standard
meridian disk of $V$ which intersects $\lambda$ in a single point. This finishes the proof of the Claim.

\QED

\noindent The annulus $ A_1 = A' \cup A''$ has the following properties: \item{\rm(i)} Both boundary 
components are on $\partial V_1 = \partial V_2$ . 
 \item{\rm(ii)} The interior of $A_1$ is in  $V_2$ . 
 \item{\rm(iii)}There is an essential  disk in $V_1$ intersecting $A_1$ in a single point.

Hence $A_1$ satisfies all the conditions for Proposition \tunnelrepProp\/ and we can change the
Heegaard splitting $(V_1, V_2)$ by replacing the tunnel $t_2$ by a tunnel $t_3$ which is an essential
arc of $A_1$. Notice that we can do this so that the new tunnel is slightly pushed off $A_1$ so it
is disjoint from it (see Fig. 4). We now have a new Heegaard splitting $(V'_1, V'_2)$ for $E(K)$ and 
the decomposing annulus $A$ intersects $V'_1$ in exactly two disks as the only change took place 
in $E(K_2)$ away from $A$. 

\bigskip

The annulus $A$ cannot minimize the intersection with $V'_1$ as $V'_1 - A$ has
exactly two components none of which has genus bigger than one. This contradicts Proposition
\innertunnelProp\/. If $|V'_1 \cap A| =1$  then we are in Case 3   treated below. If 
$|V'_1 \cap A| = 0$ then as in Case 1 the annulus $A$ is primitive in one of the Heegaard splittings  
induced on $E(K_1)$ or $E(K_2)$ by cutting $(V'_1, V'_2)$ along $A$ (as in Section 1).

\bigskip
\noindent $ \underline{Case \hskip3pt 3 }$:  In this case $A \cap V_2$ is a once punctured annulus.
Since $V_1 \cap E(K_1)$ has only one component of genus one all outermost disks of $\cal E$ 
must be contained in $E(K_2)$. As before, consider an outermost disk $\Delta$ in some meridian 
disk $E$  of  $V_2$. We have $\partial \Delta = \alpha \cup \beta$ where $\alpha$ is an arc of 
type III on $A$ and $\beta$ is an arc  on $\partial V_1 - A$.

Notice that in our situation $(V_1^1,V_2^1)$, the induced Heegaard splitting, is a genus two Heegaard 
splitting of $E(K_1)$ and $(V_1^2,V_2^2)$ is a genus three Heegaard splitting of $E(K_2)$.

Now do a boundary compression of $A$ along the disk $\Delta$ in $E(K_2)$ (An isotopy of
type A in the terminology of [Ja]). 
Denote  by $\tilde A$ the annulus obtained from $A$ after this isotopy. This isotopy does not
change the Heegaard splitting $(V_1, V_2)$ but does change the induced Heegaard splittings
on $E(K_i),  i = 1,2$. The boundary 
compression removes a regular neighborhood $N(\Delta)$ from $V^2_2$ and adds it to $V^1_2$
along the arc $\alpha$.  We obtain a new handlebody $V'^1_2$ of genus two in $E(K_1)$.
Since we have not changed the genus of the handlebody $V^1_2$ and cutting along $\tilde A$
does induce a Heegaard splitting on $E(K_1)$ we obtain a genus two Heegaard splitting 
of $E(K_1)$ by taking  $V'^1_1 =  E(K_1) - int V'^1_2$.

The situation on  $E(K_2)$ is slightly more complicated. The disk $\Delta \subset V^2_2$
is an essential disk. If it was not essential then $\Delta$ together with a disk on $\partial V_2$
bound a $3$-ball $B_0$. The boundary of $B_0$ is $\Delta$, a disk in $A$ and a disk on 
$\partial V_2$. The $3$-ball $B_0$ allows us to isotope $A$ off of $\alpha$ thus reducing the 
number of intersections of ${\cal E} \cap A$ in contradiction to the choice of $\cal E$ (see Fig. 9).
 
\bigskip
\vskip30pt
\hskip30pt \BoxedEPSF{Fig9.eps scaled 590}
%\hskip30pt { \picture 4.65in by 1.70in  (Fig.9 scaled 590)}

\bigskip
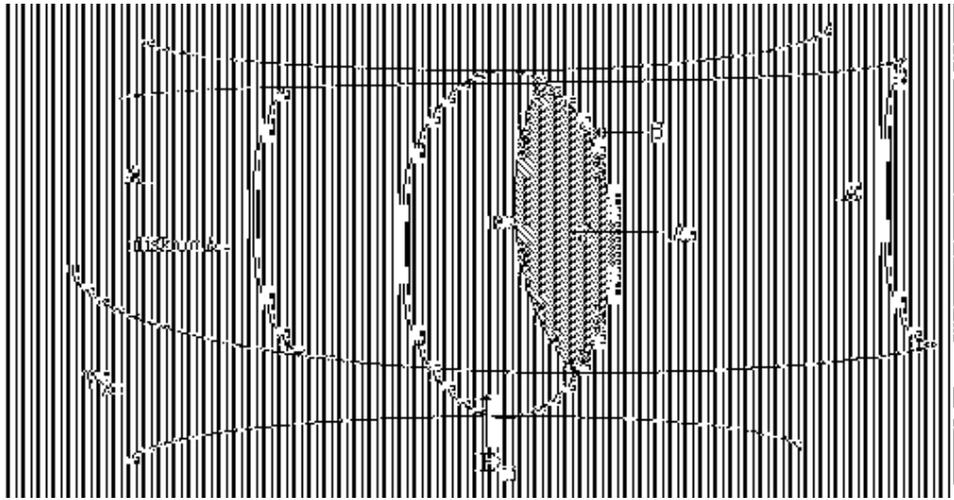
\centerline{Fig. 9} 
\bigskip
\eject

Since as before we change $A$ only by an isotopy and cutting along $\tilde A$ induces a Heegaard
splitting on $E(K_2)$ the disk $\Delta$ is non-separating and after cutting along $\tilde A$ we obtain 
a handlebody $V'^2_2$ of genus two in $E(K_2)$. Hence the  Heegaard splitting $(V'^2_1,V'^2_2)$ 
induced on $E(K_2)$ by cutting along $\tilde A$ is of genus two. 

 The annulus  $\tilde A$ intersects the
Heegaard splitting  $(V_1,V_2)$ as follows: $ \tilde A \cap V_1$ is two vertical annuli $A^*_1$
and $A^*_2$ and one essential annulus, $ \tilde A \cap V_2$ is two essential annuli 
denoted  by $A_1$ and $A_2$. Both of $A_1$ and $A_2$ are meridional annuli (see Fig. 10).
Note that the isotopy of $A$ can be done in a small neighborhood of the arc $\alpha$ on $A$
disjoint from the core curves of $A_1$ and $A_2$ .
Hence on each of the handlebodies $V'^1_2,V'^2_2$, in the induced Heegaard splittings
of $E(K_1)$ and $E(K_2)$, we have copies of the meridional annuli $A_1$ and $A_2$.

\bigskip
\bigskip
\hskip100pt \BoxedEPSF{Fig10.eps scaled500}
%\hskip100pt { \picture 2.0in by 1.80in  (Fig.10 scaled 500)}
\bigskip
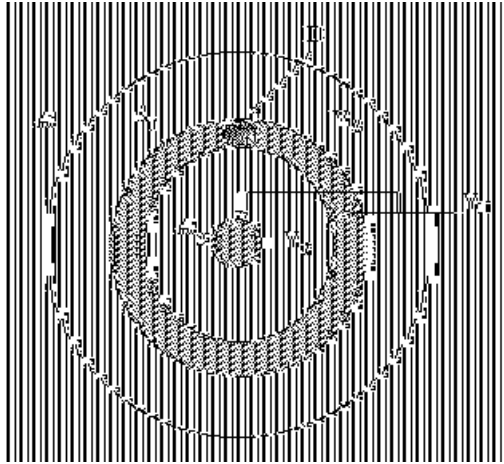
\centerline {Fig. 10}
\bigskip
Assume now that no genus two Heegaard splitting of $E(K_1)$ has a primitive meridian curve. 
Since gluing $V'^1_2$ to $V'^2_2$ along the annuli $A_1$ and $A_2$ yields a 
handlebody $V_2$ both annuli $A_1$ and $A_2$ must be primitive on $V'^2_2$. However
a priori they need not be primitive simultaneously.

Lemma 2.3.2 of [CGLS] states that two primitive curves on a genus two handlebody are either
simultaneously primitive or when we add two disks to the handlebody along the curves 
(i.e., $2$-surgery along the curves)  we obtain a non-trivial punctured Lens space. In our situation the
two curves in question are the core curves of the meridional annuli $A_1$ and $A_2$. Hence 
if the curves are not simultaneously primitive then $2$-surgery along them yields a Lens space.
But $K_2$ is a knot in $S^3$ and $2$-surgery along meridian curves yields $S^3$  in contradiction.

Recall that $V'^2_2$ was obtained from $V^2_2$ by cutting $V^2_2$ open along the disk $\Delta$.  
So in order to obtain $V^2_2$ back we need to identify the two copies of $\Delta$ on $\partial V'^2_2$.
Since $\beta \subset \partial \Delta$ and $\beta$ is paralle to an arc, also denoted by $\beta$, in
$\partial A_i$ this identification will result in
identifying the primitive annuli $A_1$ and $A_2$ on $\partial V'^2_2$ along the arc $\beta$
(see Fig. 11). Two annuli identified along an arc yield a once punctured annulus. In our situation
this planar surface has two boundary components corresponding to simultaneously  primitive 
curvs on $V'^2_2$ and the third, which is the result of the identification, corresponds to $\partial D $
where $D$ is the single disk of $V_1 \cap A$. 
\bigskip
\vskip10pt 
\hskip40pt \BoxedEPSF{Fig11.eps scaled 590}
%\hskip40pt { \picture 4.50in by 1.70in  (Fig.11 scaled 590)}
 
\bigskip

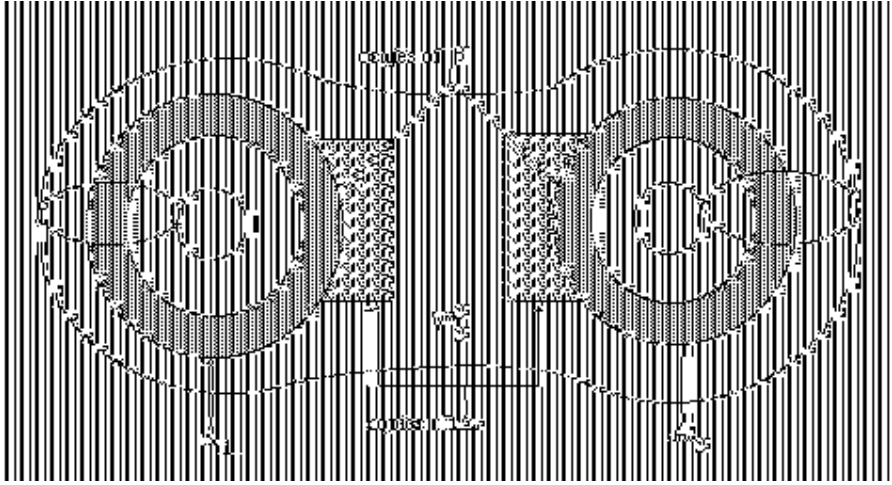
\centerline {Fig. 11}
\bigskip
Denote the generators of $\pi_1(V'^2_2)$   corresponding to the core curves of the primitive 
annuli $A_1$ and $A_2$ by $x$ and $y$ respectively.  Identifying the two copies of $\Delta \subset
\partial V'^2_2$  \lq\lq creates the third genus" of  $V_2^2$ and induces an HNN extension of
 $\pi_1(V'^2_2)$  where the new generator resulting from the identification is denoted by $z$.

Hence the third puncture in $A \cap V_2$ (i.e., $\partial D$) corresponds to the word
$x z y^{-1} z ^{-1} \in   F(x,y,z) = \pi_1(V^2_2)$. As both  $x$ and $y$ are primitive so is
$x z y^{-1} z ^{-1}$ since $\{x, z, x z y^{-1} z ^{-1} \} $ is a basis  for $F(x, y, z) $.

We conclude that  $\partial D$ is a primitive element in $\pi_1(V^2_2)$ and so must intersect
an essential disk of $V^2_2$ in a single point. Hence the Heegaard splitting of $E(K_2)$ induced
by cutting along $A$ is reducible. Thus if we remove  a regular neighborhood $N(D)$ from $V_1^2$
and add it to $V^2_2$ we still have a handlebody and a genus two Heegaard splitting of $E(K_2)$.
In the fundamental group of this handlebody $x z y^{-1} z ^{-1} = 1$ so $y = z^{-1} x z \in F(x, z)$.
But now $x =  [A_1]$ is a primitive element in this new handlebody. Hence we found a genus
two Heegaard splitting of $E(K_2)$  in which a meridional annuls is  primitive.
\smallskip
\noindent $ \underline{Case  \hskip3pt 4}$: The proof in this case is the same as for Case 3.
\bigskip
This finishes the proof of the theorem.
\QED

\bigskip
\noindent{\bf Question:} Is the statement of  Theorem  0.1 true for knots with tunnel number 
bigger than one?

\bigskip
\bigskip

\section References.

\bigskip
\noindent [CGLS] \hskip 3pt  M.Culler, C. Gordon, J. Luecke, P.Shalen; {\it Dehn surgery on knots},
Ann. of Math.,

 \noindent \hskip 39pt 125 (1987), 237 - 300

\noindent [Ja] \hskip 19pt W. Jaco; {\it Lectures on three manifold topology} CBMS Regional 
 Conference  Series  

\noindent \hskip 39pt in Mathematics 43 1977

\noindent [Kl] \hskip 19pt A. Klebanov; {\it Heegaard splittings of knot complements}, M.Sc. Thesis
Thecnion, 

\noindent \hskip 39pt Haifa, 1998

\noindent [Mo1]\hskip 15pt K. Morimoto; {\it  There are knots whose tunnel numbers go down under
connected 

\noindent \hskip 39pt sum}, Proc. Amer. Math. Soc. 123 (1995), 3527 - 3532  

\noindent [Mo2]\hskip 15pt K. Morimoto; {\it  On the additivity of tunnel number of knots}, Top. and 
App. 53

\noindent \hskip 39pt (1993), 37 - 66

\noindent [MS]\hskip 19pt K. Morimoto, M Sakuma; {\it  On unknotting tunnels for knots}, Math. Ann.
289
 
\noindent \hskip 40pt  ( 1991), 143 - 167

\noindent [MR]\hskip 19pt  Y. Moriah, H. Rubinstein; {\it Heegaard structures of negatively curved
 3-manifolds}

\noindent \hskip 40pt Comm. in Anal. and Geom. 5 (1997), 375 - 412

\noindent [MSY] \hskip 10pt K. Morimoto, M. Sakuma, Y. Yokota; {\it  Examples of 
tunnel number one knots

\noindent \hskip 42pt   which have the property that \lq\lq 1 + 1 = 3"}, Math. Proc. Camb. Phil. Soc.,

\noindent \hskip 42pt 119 (1996) 113 - 118

\noindent [No]\hskip 22pt F. Norwood; {\it Every two generator knots is prime}, Proc. Amer. Math.
Soc. 86 

\noindent \hskip 39pt (1982), 143 - 147

\noindent [Sc]  \hskip 22pt J. Schultens; {\it  Additivity of tunnel number for small knots}, preprint.

\eject

\bigskip

\obeyspaces   Yoav Moriah

\obeyspaces   Department of Mathematics

 \obeyspaces  Technion, Haifa  32000,

\obeyspaces  Israel

 ymoriah@tx.technion.ac.il

\end

%%%%%%%%%%%%%%%%%%%%%%%%%%%%%%% END %%%%%%%%%%%%%%

\bigskip
\bigskip
\bigskip

and all have the same orientaion relative to $A$ 
then $\partial \Dleta$ is a $(p,q)$ curve on the solid torus in $S^3$ bounded by the 
(decomposing $2$-sphere) $\cup$ $2$-float  giving rize to a Lens space summand of $S^3$ which
is a contradiction. 
\bigskip
%\noindent[He]  \hskip 17pt J. Hemple; {\it 3 - manifolds}, Ann. Math. Studies 86, (1976)
%Princeton University

%\hskip 26pt Press